\documentclass{article}
\usepackage{amsmath, amsfonts, amssymb}

\def\normamedida#1{ \rule[-0.8mm]{0.4mm}{3.4mm} \, #1 \,  \rule[-0.8mm]{0.4mm}{3.4mm} \:  }

\newtheorem{theorem}{Theorem}[section]
\newtheorem{lemma}[theorem]{Lemma}
\newtheorem{proposition}[theorem]{Proposition}
\newtheorem{corollary}[theorem]{Corollary}
\newtheorem{remark}[theorem]{Remark}

\numberwithin{equation}{section}

\def\dis{\displaystyle}
\newcommand{\proof}{\noindent \textit{Proof. }}
\newcommand{\qed}{\hfill $\Box$ \vspace{0.2cm}}
\newcommand{\vdois}{\vspace{0.2cm}}
\newcommand{\ep}{\varepsilon}

\begin{document}

\author{
 \textbf{Marcelo F. Furtado} \\
{\small \emph{Universidade de Bras\'{\i}lia, Departamento de Matem\'atica}}\\
{\small \emph{70910-900, Bras\'{\i}lia - DF, Brazil}}\\
{\small \emph{e-mail: mfurtado@unb.br}} \\ \\
 \textbf{Jo\~ao Pablo P. Silva} \\
{\small \emph{Universidade Federal do Par\'a, Faculdade de Matem\'atica}}\\
{\small \emph{66075-110, Bel\'em - PA, Brazil}}\\
{\small \emph{jpablo\_ufpa@yahoo.com.br}} \\ \\
}
\title{Multiplicity of solutions for homogeneous elliptic systems with critical growth \thanks{M.F. Furtado was partially supported by CNPq/Brazil
305984/2007-1.}}

\date{}

\maketitle{}

\begin{abstract}

In this paper we are concerned with the number of nonnegative solutions of the elliptic system
$$
\left\{
\begin{array}{ll}
-\Delta u = Q_u(u,v) + \frac{1}{2^*} H_u(u,v),& \mbox{ in } \Omega,\vdois\\
-\Delta v = Q_v(u,v) + \frac{1}{2^*} H_v(u,v),& \mbox{ in } \Omega,\vdois\\
u=v=0,& \mbox{on } \partial\Omega,
\end{array}
\right.
\leqno{(P)}
$$
where $\Omega \subset \mathbb{R}^N$ is a bounded smooth domain, $N \geq 3$, $2^*:=2N/(N-2)$ and $Q_u, H_u$ and $Q_v$, $H_v$ are the partial derivatives of the homogeneous functions $Q,\,H \in C^1(\mathbb{R}^2_+,\mathbb{R})$, where $
\mathbb{R}^2_+ := [0,\infty) \times [0,\infty)$.
In the proofs we apply variational methods and Ljusternik-Schnirelmann theory.

\vspace{0.2cm} \noindent \emph{2000 Mathematics Subject
Classification} : 35J20, 35J50, 58E05.

\end{abstract}

\section{Introduction}

In this paper we are concerned with the number of nonnegative solutions of the elliptic system
$$
\left\{
\begin{array}{ll}
-\Delta u = Q_u(u,v) + \frac{1}{2^*} H_u(u,v),& \mbox{ in } \Omega,\vdois\\
-\Delta v = Q_v(u,v) + \frac{1}{2^*} H_v(u,v),& \mbox{ in } \Omega,\vdois\\
u=v=0,& \mbox{on } \partial\Omega,
\end{array}
\right.
\leqno{(P)}
$$
where $\Omega \subset \mathbb{R}^N$ is a bounded smooth domain, $N \geq 3$, $2^*:=2N/(N-2)$ and $Q_u, H_u$ and $Q_v$, $H_v$ are the partial derivatives of the homogeneous functions $Q,\,H \in C^1(\mathbb{R}^2_+,\mathbb{R})$, where $
\mathbb{R}^2_+ := [0,\infty) \times [0,\infty)$.

We are interested in the case that $H$ has critical growth. More specifically, the assumptions on $H=H(s,t)$ are the following.
\begin{itemize}
 \item[$(H_0)$] $H$ is $2^*$-homogeneous, that is,
$$
 H( \theta s, \theta t)=\theta^{2^*} H(s,t)~~~ \text{for each } \theta>0,~(s,t) \in \mathbb{R}^2_+;
$$

\item[$(H_1)$] $H_s(0,1)=0,\,H_t(1,0)=0$;

\item[$(H_2)$] $H(s,t)>0$ for each $s,t>0$;

\item[$(H_3)$] $H_s(s,t)\geq 0,\,H_t(s,t) \geq 0$ for each $(s,t) \in \mathbb{R}^2_+$;

\item[$(H_4)$] the $1$-homogeneous function $(s,t) \mapsto H(s^{1/2^*},t^{1/2^*})$ is concave in $\mathbb{R}^2_+$;

\end{itemize}

The function $Q=Q(s,t)$ is a lower order perturbation term satisfying
\begin{itemize}
 \item[$(Q_0)$] $Q$ is $q$-homogeneous for some $2 \leq q < 2^*$;

  \item[$(Q_1)$] $Q_s(0,1)=0$, $Q_t(1,0)=0$.

\end{itemize}
In order to present our results we  introduce the following numbers
\begin{equation}
\mu:=\min\left\{ Q(s,t) \,:\,s^q+t^q=1,\,s,\,t\geq 0\right\}
\label{defmu}
\end{equation}
and
\begin{equation}
\lambda:=\max\left\{ Q(s,t) \,:\,s^q+t^q=1,\,s,\,t\geq 0\right\}.
\label{deflambda}
\end{equation}

We say that a weak solution $z=(u,v) \in H_0^1(\Omega) \times H_0^1(\Omega)$  of problem $(P)$ is nonnegative if  $u,\,v\geq 0$ in $\Omega$. If $Y$ is a closed set of a topological space $Z$, we denote by
cat$_Z(Y)$ the Ljusternik-Schnirelmann category of $Y$ in $Z$,
namely the least number of closed and contractible sets in $Z$
which cover $Y$. We are now ready to state our first result.

\begin{theorem} \label{th1}
Suppose that $H$ satisfies $(H_0)-(H_4)$ and $Q$ satisfies $(Q_0)-(Q_1)$.
Then there exists $\Lambda>0$ such that the problem $(P)$ has at least \emph{cat}$_{\Omega}(\Omega)$ nonzero nonnegative solutions provided $\lambda,\,\mu \in (0,\Lambda)$.
\end{theorem}

In the proof we apply variational methods, Ljusternik-Schnirelmann theory and a technique introduced by Benci and Cerami \cite{BenCer}. It consists in making precise comparisons between the category of some sublevel sets of the associated functional with
the category of the set $\Omega$.
 In order to overcame the lack of compactness due to the critical growth of $H$ we use the ideas of
Brezis and Nirenberg \cite{BreNir}, besides the paper of Morais Filho and Souto
\cite{FilSou}, where it is proved that the number
\begin{equation} \label{def_sh}
S_{H} := \inf \left\{ \displaystyle\int\limits_{\mathbb{R}^N}
(|\nabla u|^2+ |\nabla v|^2)\textrm{d}x\,:\,u,v \in
H^1(\mathbb{R}^N),\,\displaystyle\int\limits_{\mathbb{R}^N}
H(u^+,v^+)\textrm{d}x =1 \right\}
\end{equation}
plays an important role when dealing with critical systems like $(P)$.
Actually, we use the above constant and adapt some calculations
performed in Myiagaki \cite{Miy} to localize the energy levels where the
Palais-Smale condition fails. We would like to mention that, as a byproduct of our arguments, we extend the existence
result of \cite{FilSou} for any subcritical degree of homogeneity of the perturbation $Q$ (see Theorems \ref{th_existencia} and \ref{th_existencia2}).

Notice that condition $(Q_1)$ discard examples like $Q(s,t) = s^q + t^q + st^{q-1}$ since, in this case, $Q_s(0,1) = 1$. However, we can also consider this situation if the subcritical perturbation satisfies $q>2$. More specifically, the following holds
\begin{theorem} \label{th2}
Suppose that $H$ satisfies $(H_0)-(H_4)$, $Q$ satisfies $(Q_0)$ with $q>2$ and
\begin{itemize}
\item[($\widehat{Q_1})$] $Q_s(0,1)>0$ and $Q_t(1,0)>0$.
\end{itemize}
Then there exists $\Lambda>0$ such that the problem $(P)$ has at least \emph{cat}$_{\Omega}(\Omega)$ nonzero nonnegative solutions provided $\lambda,\,\mu \in (0,\Lambda)$.
\end{theorem}

The difference when dealing with $(Q_1)$ or $(\widehat{Q_1})$ is just in the way we extend the function $Q$ to the whole $\mathbb{R}^2$.  Since we want to apply minimax methods this extension needs to be made in a smooth way. We refer to the beginning of the next section for more details about the possible extensions.

Concerning the class of nonlinearities we are considering, we present in Section \ref{secao_final} some examples of functions satisfying our hypothesis. There, we also make some comments about the possibility of proving that the solutions are positive in $\Omega$ and we state other settings in which our results hold, including the possibility of having a sum of subcritical terms with different degrees of homogeneity. As a final remark, we would like to mention that the theorems remain valid for $N=3$ if the degree of homogeneity of $Q$ satisfies $4<q<6$ (see Remark \ref{n3}).

The starting point on the study of the system $(P)$ is its scalar version
\begin{equation} \label{bnp}
-\Delta u = \theta |u|^{q-2}u + |u|^{2^*-2}u~\mbox{ in } \Omega,\,\,\,u \in H_0^1(\Omega),
\end{equation}
with $2\le q<2^*$. In a pioneer work  Brezis and Nirenberg \cite{BreNir} showed that, for $q=2$,
the existence of positive solutions is related
with the interaction between the parameter $\theta$ with the first
eigenvalue $\theta_1(\Omega)$ of the operator $(-\Delta,H_0^1(\Omega))$.
Among other results they showed that, if $q=2$, the problem has at least one positive solution provided $N \geq 4$
and $0<\theta<\theta_1(\Omega)$. They also obtained some results for the case $2<q<2^*$.

After the paper of Brezis and Nirenberg, a lot of works
dealing with critical nonlinearities have been appeared. Concerning the question of multiplicity, we recall that Rey \cite{Rey} and Lazzo \cite{Laz}
proved that, for $q=2$, the problem (\ref{bnp}) has at least cat$_{\Omega}(\Omega)$ positive solutions (see also the well known paper of
Benci and Cerami \cite{BenCer} where the subcritical case was
considered) provided $\theta>0$ is small. This result was extended for the $p$-Laplacian operator and $p\leq q<p^*$ by Alves and Ding \cite{AlvDin}. The results presented here can be viewed as versions of the papers \cite{Rey,Laz,AlvDin} for the case of systems.

As far we know, the first results for homogeneous system like $(P)$ are due to Morais Filho and Souto \cite{FilSou} (see also \cite{AlvFilSou}). After this work many results have been appeared (see \cite{Han2,Han3,Han,Tan,Jap2,Hsu,Jap1}
 and references therein). Among then, the most related with our paper if the work of Han \cite{Han}, where the author considered the case $Q(s,t)=\alpha_1 s^2 + \alpha_2 t^2$ and $H(s,t) = s^{\alpha}t^{\beta}$ with $\alpha+\beta=2^*$. His results was complemented by Ishiwata in \cite{Jap1,Jap2}, with different classes of homogeneous nonlinearities being considered. Our paper extends and/or complements the results found in
\cite{FilSou,AlvFilSou,Han,Jap1,Jap2}. Although there are some multiplicity results for systems like $(P)$ via Ljusternik-Schnirelmann theory, we do not know any article that relates the topology of $\Omega$ with the number of solutions and contains a general class of nonlinearities such as those considered here.

The paper is organized as follows. In Section 2 we present the abstract
framework of the problem, we prove a local compactness result and obtain the existence of one nonnegative solution for $(P)$.
Section 3 is devoted to the proof of some technical results concerned the properties of sequences which minimize $S_H$ and the asymptotic behavior of the minimax levels associated to the problem. Theorems \ref{th1} and \ref{th2} are proved in Section 4 and we devote the last section for some further remarks about examples and possible extensions of the results.

\section{The PS condition and an existence result}

We start this section fixing some notation. We denote $B_R(0) := \{ x \in \mathbb{R}^N : \|x\| < R\}$ and by
$C_0^{\infty}(A)$ the set of all functions $f: A \to \mathbb{R}$
of class $C^{\infty}$ with compact support contained in the open
set $A \subset \mathbb{R}^N$. We denote by $\|f\|_{p}$  the $L^p$-norm of $f \in L^p(A)$.  In order to simplify the
notation, we write $\int_{A} f$ instead of $\int_{A}
f(x)\textrm{d}x$. We also omit the set $A$ whenever
$A=\Omega$.

We remark for future reference that, if $p \geq 1$ and $F$ is a $p$-homogeneous $C^1$-function, then the following holds
\begin{itemize}
\item[(i)] if we set $M_F := \max\{ F(s,t) : s,\, t \in \mathbb{R}, \,|s|^p+|t|^p=1\}$ then, for each $(s,t) \in \mathbb{R}^2$, we have that
\begin{equation} \label{prop-homogenea1}
|F(s,t)| \leq M_F(|s|^p+|t|^p)\,;
\end{equation}

\item[(ii)] $\nabla F$ is a $(p-1)$-homogeneous function and, for each $(s,t) \in \mathbb{R}^2$, we have that
\begin{equation} \label{prop-homogenea2}
sF_s(s,t) + tF_t(s,t) = pF(s,t).
\end{equation}
\end{itemize}

We proceed now with the extension of the functions $Q$ and $H$. Notice that $(Q_1)$ and $(H_1)$ allow us to give a $C^1$ extension of $Q$ and $H$ to the whole $\mathbb{R}^2$ as
\begin{equation} \label{extensao1}
\widetilde{Q}(s,t):=Q(s^+,t^+),\,\,\, \widetilde{H}(s,t):=H(s^+,t^+),
 \end{equation}
where $s^+ := \max\{s,0\}$. In the setting of Theorem \ref{th2}, with $Q$ satisfying $(\widehat{Q_1})$ instead of $(Q_1)$, the above extension is not differentiable. So, in this case, we extend $Q$ in the following way
\begin{equation} \label{extensao2}
\widetilde{Q}(s,t) := Q(s^+,t^+) - \nabla Q(s^+,t^+)\cdot(s^-,t^-),
 \end{equation}
where $s^-=\max\{-s,0\}$. For simplicity, we shall write only $H$ to denote the extension $\widetilde{H}$. The extension of $Q$ depends on it to satisfy $(Q_1)$ or $(\widehat{Q_1})$. In both cases, the extension is of class $C^1$ and will be denoted just by $Q$.

By using (\ref{prop-homogenea1}) and well know arguments, we see that the weak solutions of $(P)$ are
precisely the critical points of the $C^1$-functional
$I_{\lambda,\mu}:X \to \mathbb{R}$ given by
$$
I_{\lambda,\mu}(z) := \frac{1}{2} \|z\|^2 - \int Q_{\lambda,\mu}(z) -
\frac{1}{2^*} \int H(z),\,\,z \in X,
$$
where $X$ is the Sobolev space
$H_0^1(\Omega) \times H_0^1(\Omega)$ endowed with the norm
$$
\|(u,v)\|^2 := \int \left( |\nabla u|^2 + |\nabla v|^2 \right).
$$
We notice that, in the definition of $I_{\lambda,\mu}$, we are denoting $Q_{\lambda,\mu}(z):=Q(z)$ for $z \in \mathbb{R}^2$. We shall write $Q_{\lambda,\mu}$ instead of $Q$ just to emphasize that the smallness condition in the statement of the main theorems depends on the value of the parameters $\mu$ and $\lambda$ defined in (\ref{defmu})-(\ref{deflambda}).

We introduce the Nehari manifold of $I_{\lambda,\mu}$ by setting
$$
\mathcal{N}_{\lambda,\mu} := \left\{ z \in X \setminus \{(0,0)\} :
 I'_{\lambda,\mu}(z)z =0\right\}
$$
and define the minimax $c_{\lambda,\mu}$ as
$$
c_{\lambda,\mu} := \inf_{z\in\mathcal{N}_{\lambda,\mu}}
I_{\lambda,\mu}(z).
$$

In what follows, we present some properties of $c_{\lambda,\mu}$
and $\mathcal{N}_{\lambda,\mu}$. Its proofs can be done as in
\cite[Chapter 4]{Wil}. First of all, we note that there exists
$r=r_{\lambda,\mu}>0$, such that
\begin{equation}
\Vert z \Vert \geq r > 0~~\text{ for each } z \in
\mathcal{N}_{\lambda,\mu}.
 \label{nehari_naozero}
\end{equation}
It is standard to check that $I_{\lambda,\mu}$ satisfies Mountain
Pass geometry. So, we can use the homogeneity of $Q$ and $H$ to
prove that $c_{\lambda,\mu}$ can be alternatively characterized by
\begin{equation}
c_{\lambda,\mu} = \inf_{\gamma \in \Gamma_{\lambda,\mu}} \max_{t
\in [0,1]} I_{\lambda,\mu}(\gamma(t)) = \inf_{z \in X \setminus\{
0 \}} \max_{t \geq 0} I_{\lambda,\mu}(tz)
> 0,
 \label{prop_nehari}
\end{equation}
where $\Gamma_{\lambda,\mu} := \{ \gamma \in C([0,1],X) :
\gamma(0)=0, ~I_{\lambda,\mu}(\gamma(1)) < 0\}$. Moreover, for each
$z \in X \setminus \{0\}$, there exists a unique $t_z>0$
such that $t_z z \in \mathcal{N}_{\lambda,\mu}$. The maximum of the
function $t \mapsto I_{\lambda,\mu}(tz)$, for $t\geq 0$, is achieved
at $t=t_z$.

Let $E$ be a Banach space and $J \in C^1(E,\mathbb{R})$. We say
that $(z_n) \subset E$ is a Palais-Smale sequence at level $c$
((PS)$_c$ sequence for short) if $J(z_n) \to c$ and $J'(z_n) \to
0$. We say that $J$ satisfies (PS)$_c$ if any (PS)$_c$ sequence
possesses a convergent subsequence.

\begin{lemma}
The functional $I_{\lambda,\mu}$ satisfies the $(PS)_c$ condition
for all $c < \frac{1}{N}S_{H}^{N/2}$.
\label{lemaPS}
\end{lemma}

\proof Let $(z_n)=((u_n,v_n)) \subset X$ be such that
$I_{\lambda,\mu}'(z_n) \to 0$ and $I_{\lambda,\mu}(z_n) \to
c<\frac{1}{N}S_H^{N/2}$. The definition of $I_{\lambda,\mu}$ and (\ref{prop-homogenea1}) provide $c_1,\,c_2>0$ such that
\begin{equation} \label{psbound}
\begin{array}{lcl}
c+c_1\| z_n\| + o_n(1) & = & I_{\lambda,\mu}(z_n) - \dis\frac{1}{2^*} I_{\lambda,\mu}'(z_n)z_n \vdois \\
&=& \left( \dis\frac{1}{2} - \dis\frac{1}{2^*} \right) \|z_n\|^2 + \left( \dis\frac{q-2^*}{2^*} \right) \dis\int Q_{\lambda,\mu}(z_n) \vdois \\
& \leq & c_2 \left( \|z_n\|^2 + \|z_n\|^q\right),
\end{array}
\end{equation}
where hereafter $o_n(1)$ denotes a quantity approaching zero as
$n\to\infty$. The above expression implies that $(z_n) \subset X$
is bounded. So, we may suppose that $z_n \rightharpoonup z:=(u,v)$
weakly in $X$ and $z_n \to z$ strongly in $L^q(\Omega) \times
L^q(\Omega)$. Moreover, a standard argument shows that
$I_{\lambda,\mu}'(z)=0$.

By setting $\widetilde{z}_n := (\widetilde{u}_n,\widetilde{v}_n) =
(u_n-u,v_n-v)$ we can use the strong convergence in $L^q(\Omega)
\times L^q(\Omega)$ and \cite[Lemma 5]{FilSou} to conclude that
\begin{equation}
\int Q_{\lambda,\mu}(z_n) = \int Q_{\lambda,\mu}(z) + o_n(1),~~~\int H(z_n) =\int H(z) + \int
H(\widetilde{z}_n)+o_n(1). \label{brezislieb}
\end{equation}
This and the weak convergence of $(z_n)$ provide
\begin{equation} \label{ps1}
c+o_n(1) = I_{\lambda,\mu}(z) + \frac{1}{2} \| \widetilde{z}_n\|^2
- \frac{1}{2^*} \int H(\widetilde{z}_n) \geq \frac{1}{2} \|
\widetilde{z}_n\|^2 - \frac{1}{2^*} \int H(\widetilde{z}_n),
\end{equation}
where we have used $I_{\lambda,\mu}(z) \geq 0$.

By using $I_{\lambda,\mu}'(z_n) \to 0$ and (\ref{brezislieb}) again, we get
$$
\begin{array}{lcl}
o_n(1) &=& I_{\lambda,\mu}'(z_n)z_n = \| z_n\|^2 - q \dis\int Q_{\lambda,\mu}(z_n) - \dis\int H(z_n) \vdois \\
 & = & I_{\lambda,\mu}'(z)z + \| \widetilde{z}_n\|^2 - \dis\int H(\widetilde{z}_n).
\end{array}
$$
Recalling that $I_{\lambda,\mu}'(z)=0$, we can use the above equality and (\ref{ps1}) to obtain
$$
\lim_{n\to \infty} \| \widetilde{z}_n\|^2 = b = \lim_{n\ \to
\infty} \int H(\widetilde{z}_n),~~~\frac{1}{N} b = \left(
\frac{1}{2}-\frac{1}{2^*} \right)b \leq c,
$$
for some $b \geq 0$.

In view of the definition of $S_H$, we have that
$$
\| \widetilde{z}_n \|^2 \geq S_H \left( \int H(\widetilde{z}_n ) \right)^{2/2^*}.
$$
Taking the limit we get $b \geq S_H b^{2/2^*}$. So, if $b >0$, we conclude that $b \geq S_H^{N/2}$ and therefore
$$
\frac{1}{N} S_H^{N/2} \leq \frac{1}{N} b \leq c < \frac{1}{N} S_H^{N/2},
$$
which does not make sense. Hence $b=0$ and therefore $z_n \to z$ strongly in $X$. \qed

Before presenting our next result we recall that, for each $\ep>0$, the function
\begin{equation}
\Phi_{\ep}(x):=  \frac{C_N \ep^{(N-2)/4}}{(\ep+|x|^2)^{(N-2)/2}},~~~x
\in\mathbb{R}^N, \label{instanton}
\end{equation}
where $C_N := N(N-2)^{(N-2)/4}$, satisfies $\|\nabla
\Phi{_\ep}\|_2^2=\|\Phi_{\ep}\|^{2^{*}}_{2^*}=S^{N/2}$, where $S$ is the best constant of the Sobolev embedding $\mathcal{D}^{1,2}(\mathbb{R}^N) \hookrightarrow L^{2^*}(\mathbb{R}^N)$. Thus,
using \cite[Lemma 1]{FilSou} and the homogeneity of $H$,  we
obtain $A,\,B>0$ such that
$$
S_H  = \dis\frac{||(A \Phi_{\ep},B \Phi_{\ep})||^2}{\left(\dis\int_{\mathbb{R}^N} H(A \Phi_{\ep},B \Phi_{\ep})\right)^{2/2{*}}} = \dis\frac{(A^2+B^2)}{H(A,B)^{2/2^{*}}}
\dis\frac{S^{N/2}}{\|\Phi_{\ep}\|_{2^*}^2},
$$
from which it follows that
\begin{equation}
S_H =
\frac{(A^2+B^2)}{H(A,B)^{2/2^{*}}}S. \label{SH}
\end{equation}

The above equality and the ideas introduced by Brezis and Nirenberg \cite{BreNir} are the keystone of the following result.

\begin{lemma}
Suppose that $Q$ satisfies $(Q_0)$, with $2<q<2^*$, and $\lambda,\,\mu$ defined in (\ref{defmu})-(\ref{deflambda}) are positive. Then,
$$
c_{\lambda,\mu} <
\frac{1}{N}S_{H}^{N/2}.
$$
The same result holds if $q=2$ and
and $\lambda,\,\mu \in (0,\theta_1(\Omega)/2)$, where $\theta_1(\Omega)>0$
denotes the first eigenvalue of $(-\Delta,H_0^1(\Omega))$.
\label{lema-minimax}
\end{lemma}

\proof We consider a nonnegative function $\phi \in
C_0^{\infty}(\mathbb{R}^N)$ such that $\phi \equiv 1$ in $B_R(0)
\subset \Omega$, $\phi \equiv 0$ in $\mathbb{R}^N \setminus
B_{2R}(0)$ and define
$$
w_{\ep}(x) := \frac{\phi(x) \Phi_{\ep}(x)}{\| \phi,
\Phi_{\ep}\|_{2^*}}.
$$
where $\Phi_{\ep}$ was defined in (\ref{instanton}).
Since $\| w_{\ep} \|_{2^*}=1$, we can use the homogeneity of  $Q$ and $H$ to get, for any $t \geq 0$,
$$
I_{\lambda,\mu}(tAw_{\ep},tBw_{\ep}) = \frac{t^2}{2}(A^2+B^2) \|w_{\ep}\|^2 - t^qQ_{\lambda,\mu}(A,B) \| w_{\ep}\|^q_q - \frac{t^{2^*}}{2^*}H(A,B).
$$
We shall denote by $h_{\ep}(t)$ the right-hand side of the above equality and consider two distinct cases.

\vdois\noindent \textbf{Case 1}. $2<q<2^*$.
\vdois

\noindent In this case there exists $t_{\ep}>0$ such that
\begin{equation} \label{ps4}
h_{\ep}(t_{\ep})= \max_{t \geq 0} h_{\ep}(t).
\end{equation}
Let
$$
g_{\ep}(t):=\dis\frac{t^2}{2}(A^2+B^2)\|w_{\ep}\|^2-\dis\frac{t^{2^{*}}}{2^{*}} H(A, B),~~~ t \geq 0,
$$
and notice that the maximum value of $g_{\ep}$ occurs at the point
$$
\tilde{t}_{\ep} := \left\{\frac{(A^2+B^2)\|w_\ep\|^2}{H(A,B)}\right\}^{1/(2-2^{*})}.
$$
So, for each $t\geq 0$,
$$
g_{\ep}(t)\leq g_{\ep}(\tilde{t}_{\ep}) = \frac{1}{N}
\left(\frac{(A^2+B^2)\|w_\ep\|^2}{H(A,B)^{2/2^{*}}}\right)^{N/2},
$$
and therefore
\begin{equation} \label{ps2}
h_{\ep}(t_\ep )\leq \frac{1}{N}
\left(\frac{(A^2+B^2)\|w_\ep\|^2}{H(A,B)^{2/2^{*}}}\right)^{N/2}-t_{\ep}^q
Q_{\lambda,\mu}(A,B)\|w_\ep\|^q_q.
\end{equation}

We claim that, for some $c_2>0$, there holds
$$
t_{\ep}^q Q_{\lambda,\mu}(A,B) \geq c_2.
$$
Indeed, if this is not the case, we have that $t_{\ep_n} \to 0$ for some sequence $\ep_n \to 0^+$. But it is proved in \cite[(1.11) and (1.12)]{BreNir}
that
\begin{equation}
\|w_\ep\|^2=S+O(\ep^{(N-2)/2}).
\label{S}
\end{equation}
Thus,
$$
0< c_{\lambda,\mu} \leq \sup_{t \geq 0}
I_{\lambda,\mu}(tAw_{\ep_n},tBw_{\ep_n}) =
I_{\lambda,\mu}(t_{\ep_n}Aw_{\ep_n},t_{\ep_n}Bw_{\ep_n}) \to 0,
$$
which does not make sense. So, the claim holds and
we infer from (\ref{ps2}) and (\ref{S}) that
$$
\begin{array}{lcl}
h_\ep(t_\ep)&\leq&
\dis\frac{1}{N}\left(\dis\frac{(A^2+B^2)}{H(A,B)^{2/2^{*}}}
S+O(\ep^{(N-2)/2})\right)^{N/2}-c_2 \|w_\ep\|^q_q
\vspace{0.3cm}
\\
&\leq&\dis\frac{1}{N}S^{N/2}_H+O(\ep^{(N-2)/2})-c_2 \|w_\ep\|^q_q.
\end{array}
$$
It is proved in \cite[Claim 2, p. 778]{Miy} that $\lim_{\ep \to 0^+} \ep^{(2-N)/2}\|w_{\ep}\|_q^q = + \infty$. Thus, we conclude from the above inequality that, for each $\ep>0$ small, there holds
$$
c_{\lambda,\mu} \leq \sup_{t \geq 0}
I_{\lambda,\mu}(tAw_{\ep},BAw_{\ep}) = h_{\ep}(t_{\ep}) <
\frac{1}{N} S_H^{N/2}.
$$

\vdois\noindent \textbf{Case 2}. $q=2$.
\vdois

\noindent In this case we have that $h_\ep^{'}(t)=0$ if, and only if,
$$
(A^2+B^2)\|w_\ep\|^2-2Q_{\lambda,\mu}(A,B) \|w_\ep\|^2_2=t^{2^{*}-2}H(A,B).
$$
Since we are supposing $\lambda<\theta_1(\Omega)/2$, we can use Poincar\'e's Inequality to obtain
$$
\begin{array}{lcl}
2Q_{\lambda,\mu}(A,B) \| w_{\ep}\|^2_2 &\leq& 2\lambda (A^2+B^2) \| w_{\ep}\|^2_2 \vdois \\
&<&
\theta_1(\Omega)(A^2+B^2)\| w_{\ep}\|^2_2 \leq (A^2+B^2) \| w_{\ep}\|^2.
\end{array}
$$
Thus, there exists $t_{\ep}>0$ satisfying (\ref{ps4}).
By using the definition of $w_{\ep}$ and \cite[(1.12) and (1.13)]{BreNir} we get
\begin{equation} \label{bn}
\|w_\ep\|^2_2 =
\left\{
\begin{array}{ll}
\ep^{(N-2)/4}+O(\ep^{(N-2)/2}) & \mbox{if }N\geq 5,\vdois\\
\ep^{(N-2)/2}|\log \ep|+O(\ep^{(N-2)/2}) & \mbox{if }N=4.
\end{array}
\right.
\end{equation}
Arguing as in the first case we conclude that, for $\ep>0$ small, there holds
$$
h_\ep(t_\ep)\leq\frac{1}{N}S_H^{N/2}+O(\ep^{(N-2)/2})-c_2\|w_\ep\|^2_2 < \frac{1}{N}S_H^{N/2},
$$
where we have used (\ref{bn}) in the last inequality. This concludes the proof. \qed

\begin{remark} \label{n3}
The previous lemma remains valid if we suppose that $N=3$ and $4<q<6$. Indeed, it suffices to notice that in this case, according to \cite[p. 779]{Miy}, the function $w_{\ep}$ above satisfies $\lim_{\ep \to 0^+} \ep^{(2-N)/2}\|w_{\ep}\|_q^q = + \infty$. So, the same arguments of Case 1 hold.
\end{remark}

As a byproduct of Lemmas \ref{lemaPS} and \ref{lema-minimax} we obtain the following generalization of \cite[Theorem 1]{FilSou}.

\begin{theorem} \label{th_existencia}
Under the hypotheses of Theorem \ref{th1} the problem $(P)$ possesses a nonzero nonnegative solution whenever $2<q<2^*$
and $\lambda,\,\mu>0$, or $q=2$ and $\lambda,\,\mu \in
(0,\theta_1(\Omega)/2)$. The same result holds if $N=3$ and $4<q<6$.
\end{theorem}

\proof Since $I_{\lambda,\mu}$ satisfies the geometric conditions of the Mountain Pass Theorem, there exists $(z_n) \subset X$ such that
$$
I_{\lambda,\mu}(z_n) \to c_{\lambda,\mu},~~~I_{\lambda,\mu}'(z_n)
\to 0.
$$
It follows from Lemma \ref{lemaPS} (with Remark \ref{n3} in the case $N=3$) and Lemma \ref{lema-minimax} that $(z_n)$ converges, along a
subsequence, to a nonzero critical point $z=(u,v) \in X$ of
$I_{\lambda,\mu}$. According to (\ref{extensao1}) and (\ref{prop-homogenea2}), we have that
$$
I_{\lambda,\mu}(z)z^- =  -\| z^-\|^2 - \dis\int \left(\nabla Q(u^+v^+)\cdot (u^-v^-) + \frac{1}{2^*} \nabla H(u^+v^+)\cdot (u^-v^-) \right)
$$
Since $z$ is a critical point and the integral above vanishes, it follows that $z^-=0$. Hence, $u,\,v \geq 0$ in $\Omega$ and the theorem is proved. \qed

\begin{theorem} \label{th_existencia2}
Under the hypotheses of Theorem \ref{th2} the problem $(P)$ possesses a nonzero nonnegative solution whenever  $\lambda,\,\mu>0$.  The same result holds if $N=3$ and $4<q<6$.
\end{theorem}

\proof As before, we obtain a nonzero critical point $z$ of $I_{\lambda,\mu}$. A simple calculation shows that the extension given in (\ref{extensao2}) implies that $Q_s(s,t) \geq 0$ for $s \leq 0$, and $Q_t(s,t) \geq 0$ for $t \leq 0$. Hence, using the extension of $H$ and arguing as in the previous theorem we obtain
$$
0 = I_{\lambda,\mu}'(z)z^- = -\| z^-\|^2 - \int \left( Q_u(u,v)u^- + Q_v(u,v)v^- \right)  \leq -\| z^-\|^2,
$$
and the result follows. \qed

\section{Some technical results}

In this section we denote by ${\mathcal{M}}(\mathbb{R}^N)$ the
Banach space of finite Radon measures over ${\mathbb{R}^N}$
equipped with the norm
$$
\normamedida{\sigma} = \displaystyle\sup_{\varphi \in
C_0(\mathbb{R}^N),\|\varphi\|_{\infty} \leq 1} |\sigma(\varphi)|.
$$
A sequence $(\sigma_n) \subset {\mathcal{M}}(\mathbb{R}^N)$ is
said to converge weakly to $\sigma \in
{\mathcal{M}}(\mathbb{R}^N)$ provided $\sigma_n(\varphi) \to
\sigma(\varphi)$ for all $\varphi \in C_0(\mathbb{R}^N)$. By the
Banach-Alaoglu theorem, every bounded sequence $(\sigma_n) \subset
{\mathcal{M}}(\mathbb{R}^N)$ contains a weakly convergent
subsequence.

The next result is a version of the Second
Concentration-Compactness Lemma of P.L.Lions \cite[Lemma
I.1]{Lio2}.

\begin{lemma}
Suppose that the sequence $(w_n) \subset
\mathcal{D}^{1,2}(\mathbb{R}^N) \times \mathcal
{D}^{1,2}(\mathbb{R}^N)$ satisfies
$$
\begin{array}{llc}
w_n\rightharpoonup w  & \mbox{weakly in }
\mathcal{D}^{1,2}(\mathbb{R}^N)\times
\mathcal{D}^{1,2}(\mathbb{R}^N),& \vspace{0.3cm}
\\
w_n(x) \rightarrow w(x) &\mbox{for a.e. }x \in \mathbb{R}^N,&
\vspace{0.3cm}
\\
|\nabla (w_n-w)|^2 \rightharpoonup \sigma & \mbox{weakly in }
\mathcal{M}(\mathbb{R}^N),& \vspace{0.3cm}
\\
H(w_n-w)\rightharpoonup \nu & \mbox{weakly in }
\mathcal{M}(\mathbb{R}^N)&
\end{array}
$$
and define
\begin{equation}
\sigma_{\infty}:= \lim_{R \rightarrow \infty}\limsup_{n\to\infty}
\dis\int_{|x|>R}|\nabla w_n|^2 \textrm{d}x,~~~~ \nu_{\infty} :=
\lim_{R\rightarrow \infty}\limsup_{n \to \infty}
\dis\int_{|x|>R}H(w_n)\textrm{d}x. \label{nuinfinito}
\end{equation}
Then
\begin{eqnarray}
\limsup_{n \to \infty}\dis\int_{\mathbb{R}^N} |\nabla
w_n|^2\textrm{d}x = \normamedida{\sigma}+ \sigma_{\infty}+
\dis\int_{\mathbb{R}^N} |\nabla w|^2\textrm{d}x,
\label{limsupuv} \vspace{1cm} \\
\limsup_{n \to \infty}\dis\int_{\mathbb{R}^N}
H(w_n)\textrm{d}x=\normamedida{\nu}+\nu_{\infty} +
\dis\int_{\mathbb{R}^N} H(w)\textrm{d}x,
\label{limsuphuv} \vdois \\
\normamedida{\nu}^{2/2^{*}}\leq
S_H^{-1}\normamedida{\sigma}~~\mbox{ and }~~
\nu_{\infty}^{2/2^{*}}\leq S_H^{-1}\sigma_{\infty}.
\label{desigmedidas}
\end{eqnarray}
Moreover, if $w=0$ and
$\normamedida{\nu}^{2/2^{*}}=S_H^{-1}\normamedida{\sigma}$, then
there exists $x_0,\,x_1 \in \mathbb{R}^N$ such that
$\nu=\delta_{x_0}$ and $\sigma=\delta_{x_1}$.
\label{lemmaccompacidade}
\end{lemma}

\proof We first recall that, in view of the definition of $S_H$,
for each nonnegative function $\varphi \in
C_0^{\infty}(\mathbb{R}^N)$ we have that
$$
\left( \int_{\mathbb{R}^N} \varphi^{2^*}(x)H(w_n)\textrm{d}x
\right)^{2/2^*} = \left( \int_{\mathbb{R}^N} H(\varphi(x)
w_n)\textrm{d}x \right)^{2/2^*} \leq S_H^{-1} \| \varphi(x)
w_n\|^2.
$$
Moreover, arguing as in \cite[Lemma 5]{FilSou}, we have that
$$
\int_{\mathbb{R}^N} \psi(x) H(w_n-w)\textrm{d}x =
\int_{\mathbb{R}^N} \psi(x) H(w_n)\textrm{d}x -
\int_{\mathbb{R}^N} \psi(x) H(w)\textrm{d}x + o_n(1),
$$
for each $\psi \in C_0^{\infty}(\mathbb{R}^N)$. Since $H$ is
$2^*$-homogeneous, we can use the two above expressions and argue along the same lines of the proof of \cite[Lemma 1.40]{Wil} (see
also \cite[Lemma 2.2]{Han}) to conclude that
(\ref{limsupuv})-(\ref{desigmedidas}) hold. If $w=0$ and
$\normamedida{\nu}^{2/2^{*}}=S_H^{-1}\normamedida{\sigma}\,$ the same argument of \cite[step 3 of the proof of Lemma 1.40]{Wil} implies that the measures $\nu$ and $\sigma$  are concentrated at
single points $x_0, \, x_1 \in \mathbb{R}^N$, respectively. \qed

\begin{remark} \label{remarkcompacidade}
For future reference we notice that the last conclusion of the
above result holds even  if $w \not\equiv 0$. Indeed, in this
case we can  define $\widetilde{w}_n :=
w_n-w$ and notice that
$$
\begin{array}{llc}
\widetilde{w}_n \rightharpoonup \widetilde{w} = 0  & \mbox{weakly
in } \mathcal{D}^{1,2}(\mathbb{R}^N)\times
\mathcal{D}^{1,2}(\mathbb{R}^N),& \vspace{0.3cm}
\\
\widetilde{w}_n(x) \rightarrow 0 &\mbox{for a.e. }x \in
\mathbb{R}^N,& \vspace{0.3cm}
\\
|\nabla (\widetilde{w}_n-\widetilde{w})|^2 \rightharpoonup
\widetilde{\sigma} & \mbox{weakly in } \mathcal{M}(\mathbb{R}^N),&
\vspace{0.3cm}
\\
H(\widetilde{w}_n-\widetilde{w})\rightharpoonup \widetilde{\nu} &
\mbox{weakly in } \mathcal{M}(\mathbb{R}^N).&
\end{array}
$$
But $\widetilde{w}_n - \widetilde{w} = w_n - w$ and therefore
$\widetilde{\sigma}=\sigma$ and $\widetilde{\nu} = \nu$,  where $\sigma$ and $\nu$ are as in Lemma \ref{lemmaccompacidade}. Thus, if
$\normamedida{\nu}^{2/2^{*}}=S_H^{-1}\normamedida{\sigma}$ we also
have that
$\normamedida{\widetilde{\nu}}^{2/2^{*}}=S_H^{-1}\normamedida{\widetilde{\sigma}}$
and the result follows from the last part of Lemma \ref{lemmacompacidade}.
\end{remark}

Before stating one of the main results of this section we introduce the following notation. Given
$r>0$, $y \in \mathbb{R}^N$ and a function $z \in X$, we extend
$z$ to the whole $\mathbb{R}^N$ by setting $z(x):=0$ if $x \in
\mathbb{R}^N \setminus \Omega$ and define $z^{y,r} \in
H^1(\mathbb{R}^N) \times H^1(\mathbb{R}^N)$ as
$$
z^{y,r}(x) := r^{(N-2)/2}z(rx+y),~~~x \in \mathbb{R}^N.
$$

\begin{proposition}
Suppose $(z_n) \subset X$ is such that
$$
\int H(z_n)=1~~\mbox{ and }~~\lim_{n\to\infty} \| z_n\|^2 = S_H.
$$
Then there exist $(r_n) \subset (0,\infty)$ and $(y_n) \subset
\mathbb{R}^N$ such that the  sequence $(z_n^{y_n,r_n})$
strongly converges to $z \neq 0$ in
$\mathcal{D}^{1,2}(\mathbb{R}^N) \times
\mathcal{D}^{1,2}(\mathbb{R}^N)$. Moreover, as $n\to\infty$, we
have that $r_n \to 0$ and $y_n \to \overline{y}\in
\overline{\Omega}$. \label{lemaSH}
\end{proposition}

\proof We first extend $z_n$ by setting $z_n(x):=0$ if $x\in\mathbb{R}^N \setminus \Omega$. For each $r>0$ we consider
$$
F_n(r) := \sup_{y \in \mathbb{R}^N} \int_{B_r(y)} H(z_n).
$$
Since $\lim_{r\to 0} F_n(r)=0$ and $\lim_{r\to\infty}
F_n(r)=1$, there exist $r_n>0$ and a sequence
$(y_n^k)_{k\in\mathbb{N}} \subset \mathbb{R}^N$ satisfying
$$
\frac{1}{2} =F_n(r_n) = \lim_{k\to\infty} \int_{B_{r_n}(y_n^k)}
H(z_n).
$$
Recalling that $\lim_{|y|\to\infty} \int_{B_{r_n}(y)} H(z_n)=0$ we
conclude that $(y_n^k)$ is bounded. Hence, up to a subsequence,
$\lim_{k\to\infty} y_n^k = y_n \in \mathbb{R}^N$ and we obtain
$$
\frac{1}{2} = \int_{B_{r_n}(y_n)} H(z_n).
$$

We shall prove that the sequences $(r_n)$ and $(y_n)$ above
satisfy the statements of the lemma. First notice that
\begin{equation} \label{minimiza1}
\frac{1}{2} = \int_{B_{r_n}(y_n)} H(z_n) = \int_{B_{1}(0)}
H(z_n^{y_n,r_n}) = \sup_{y \in \mathbb{R}^N}\int_{B_{1}(y)}
H(z_n^{y_n,r_n}).
\end{equation}
If we denote $w_n := z_n^{y_n,r_n}$, a straightforward calculation
provides
$$
\lim_{n\to\infty} \| w_n\|^2 = \lim_{n\to\infty} \|z_n\|^2 =
S_H,~~~\int_{\mathbb{R}^N} H(w_n)=1.
$$
Hence, we can apply Lemma \ref{lemmaccompacidade} to obtain $w \in
H^1(\mathbb{R}^N) \times H^1(\mathbb{R}^N)$ satisfying
\begin{eqnarray}
 S_H= \normamedida{\sigma}+ \sigma_{\infty}+  \|w\|^2,~~~
1=\normamedida{\nu}+\nu_{\infty} + \dis\int_{\mathbb{R}^N} H(w),
\label{minimiza2} \vdois \\
\normamedida{\nu}^{2/2^{*}}\leq
S_H^{-1}\normamedida{\sigma}~~\mbox{ and }~~
\nu_{\infty}^{2/2^{*}}\leq S_H^{-1}\sigma_{\infty}.
\label{minimiza3}
\end{eqnarray}

The second equality above implies that $\int
H(w),\,\normamedida{\nu},\, \nu_{\infty} \in [0,1]$. If one of
these values belongs to the open interval $(0,1)$, we can use
(\ref{minimiza2}), $2/2^*<1$, $(\int H(w))^{2/2^*} \leq
S_H^{-1}\|w\|^2$ and (\ref{minimiza3}) to get
$$
\begin{array}{lcl}
S_H & = & S_H \left( \normamedida{\nu}+\nu_{\infty} +
\dis\int_{\mathbb{R}^N} H(w)\right) \vdois \\
&<& S_H\left( \normamedida{\nu}^{2/2^*} + \nu_{\infty}^{2/2^*} +
\left( \dis\int_{\mathbb{R}^N} H(w) \right)^{2/2^{*}} \right) \leq
S_H,
\end{array}
$$
which does not make sense. Thus  $\int H(w),\,\normamedida{\widetilde{\nu}},\,
\nu_{\infty} \in \{0,1\}$. Actually, it follows from
(\ref{minimiza1}) that $\int_{|x|>R} H(w_n) \leq 1/2$ for any $R
>1$. Thus, we conclude that
$\nu_{\infty}=0$.

Let us prove that $\normamedida{\nu}=0$. Suppose, by
contradiction, that $\normamedida{\nu}=1$. It follows from
the first equality in (\ref{minimiza3}) that $S_H \leq
\normamedida{\sigma}$. On the other hand, the first equality in
(\ref{minimiza2}) provides $\normamedida{\sigma} \leq S_H$. Hence,
we conclude that $\normamedida{\sigma}=S_H$. Since we are
supposing that $\normamedida{\nu}=1$ we obtain $\normamedida{
\nu}^{2/2^*}=S_H^{-1}\normamedida{\sigma}$. It follows from Remark
\ref{remarkcompacidade} that $\nu=\delta_{x_0}$ for some $x_0 \in
\mathbb{R}^N$. Thus, from (\ref{minimiza1}), we get
$$
\frac{1}{2} \geq \lim_{n\to\infty}\int_{B_1(x_0)}
H(w_n) = \int_{B_1(x_0)} \textrm{d}\nu =
\normamedida{\nu}=1.
 $$
This contradiction proves that $\normamedida{\nu}=0$.

Since $\normamedida{\nu}=\nu_{\infty}=0$ we have that
$\int_{\mathbb{R}^N} H(w)=1$. This and (\ref{minimiza2}) provide
$$
\lim_{n\to\infty} \|w_n\|^2 = S_H \geq \|w\|^2 \geq
S_H\left(\int_{\mathbb{R}^N} H(w)\right)^{2/2^*}=S_H.
$$
So, $\|w\|^2=S_H$ and therefore $w_n \to w \not\equiv 0$ strongly in
$\mathcal{D}^{1,2}(\mathbb{R}^N) \times
\mathcal{D}^{1,2}(\mathbb{R}^N)$ and $w_n(x) \to w(x)$ for a.e. $x
\in\mathbb{R}^N$.

In order to conclude the proof we notice that
$$
\| w_n \|_{L^2(\mathbb{R}^N) \times L^2(\mathbb{R}^N)} =
\frac{1}{r_n^2} \| z_n \|_{L^2(\Omega) \times L^2(\Omega)}.
$$
Since $(z_n)$ is bounded and $w \not\equiv 0$, we infer from the above equality that, up to a subsequence,  $r_n \to r_0 \geq 0$. If $|y_n| \to \infty$
we have that, for each fixed $x \in \mathbb{R}^N$, there exists $n_x\in\mathbb{N}$ such that $r_nx+y_n\notin\Omega$ for $n\geq n_x$. For such values of $n$ we have that $w_n(x)=0$. Taking the limit and recalling that $x\in\mathbb{R}$ is arbitrary, we conclude that $w\equiv 0$, which is absurd. So,
along a subsequence, $y_n \to y \in \mathbb{R}^N$.

We claim that $r_0=0$. Indeed, suppose by contradiction that
$r_0>0$. Then, as $n$ becomes large, the set $\Omega_n :=
(\Omega-y_n)/r_n$ approaches to $\Omega_0 := (\Omega - y)/r_0 \neq
\mathbb{R}^N$. This implies that $w$ has compact
support in $\mathbb{R}^N$. On the other hand, since $w$ achieves the infimum in
(\ref{def_sh}) and $H$ is homogeneous, we
can use the Lagrange Multiplier Theorem to conclude that $w=(u,v)$
satisfies
$$
-\Delta u= \lambda H_u(u,v),~-\Delta v=\lambda H_v(u,v),~~x \in
\mathbb{R}^N,
$$
for $\lambda=2S_H/2^*>0$. It follows from $(H_3)$ and
the Maximum Principle that at least one the functions $u$, $v$ is positive in $\mathbb{R}^N$. But this
contradicts supp$\,w \subset \Omega_0$. Hence, we conclude that
$r_0=0$. Finally, if $y \not\in \overline{\Omega}$ we obtain $r_nx
+ y_n \not\in \Omega$ for large values of $n$, and therefore we
should have $w \equiv0$ again. Thus, $y \in \overline{\Omega}$ and the
proof is finished. \qed

We finalize this section with the study of the asymptotic behavior of the minimax level $c_{\lambda,\mu}$ as both the parameters approaches zero.

\begin{lemma} \label{lema-limite-zero}
We have that
$$
\lim_{\lambda, \,\mu \to 0^+} c_{\lambda,\mu} = c_{0,0} = \frac{1}{N}S_H^{N/2}.
$$
\end{lemma}

\proof We first prove the second equality. It follows from $\lambda=\mu=0$ that $Q_{0,0} \equiv 0$. If $A$,
$B$, $w_{\ep}$, $g_{\ep}$ and $t_{\ep}$ are as in the proof of
Lemma \ref{lemaPS}, we have that $(t_{\ep}Aw_\ep,t_{\ep}Bw_\ep)\in
\mathcal{N}_{0,0}$. Thus
$$
\begin{array}{lcl}
c_{0,0} & \leq & I_{0,0}(t_{\ep}A w_\ep,t_{\ep}B w_\ep)=\dis\frac{1}{N}\left\{ \frac{(A^2+B^2)}{H(A,B)^{2/2^{*}}} \|w_\ep\|^2\right\}^{N/2} \vdois \\
& = & \dis\frac{1}{N}\left\{ \frac{(A^2+B^2)}{H(A,B)^{2/2^{*}}} (S+O(\ep^{(N-2)/2})) \right\}^{N/2}.
\end{array}
$$
Taking the limit as $\ep \to 0^+$ and using (\ref{SH}), we conclude that $c_{0,0} \leq
\frac{1}{N}S_H^{N/2}$.

In order to obtain the reverse inequality we consider $(z_n) \subset X$ such that $I_{0,0}(z_n) \to c_{0,0}$ and $I_{0,0}'(z_n) \to 0$. The sequence $(z_n)$ is bounded and therefore
$I_{0,0}'(z_n) z_n = \| z_n\|^2 - \int H(z_n) = o_n(1)$. It follows that
$$
\lim_{n\to \infty} \| z_n\|^2 = b=\lim_{n\ \to \infty} \int
H(z_n).
$$
Taking the limit in the inequality
$
S_H \left( \int H(z_n ) \right)^{2/2^*} \leq \| z_n\|^2
$
we conclude, as in the proof of Lemma \ref{lemaPS}, that $Nc_{0,0} = b \geq S_H^{N/2}$. Hence,
$$
c_{0,0} = \lim_{n\to\infty} I_{0,0}(z_n) = \lim_{n\to \infty} \left( \frac{1}{2} \| z_n\|^2 - \frac{1}{2^*} \int H(z_n)\right) = \frac{1}{N} b \geq \frac{1}{N}S_H^{N/2},
$$
and therefore $c_{0,0}=\frac{1}{N}S_H^{N/2}$.

We proceed now with the calculation of $\lim_{\lambda,\,\mu \to 0^+} c_{\lambda,\mu}$. Let $(\lambda_n),\,(\mu_n) \subset \mathbb{R}^+$ be such that $\lambda_n,\,\mu_n \to 0^+$. Since $\mu_n$ defined in (\ref{defmu}) is positive, we have that $Q_{\lambda_n,\mu_n}(z)\geq 0$ whenever $z$ is nonnegative. Thus, for this kind of function, we have that $I_{\lambda_n,\mu_n}(z) \leq I_{0,0}(z)$. It follows that
$$
\begin{array}{lcl}
c_{\lambda_n,\mu_n}&=& \dis\inf_{z\neq (0,0)}\max_{t\geq0}I_{\lambda_n,\mu_n}(tz)
\vspace{0.2cm}
\\
&\leq& \dis\inf_{z \neq (0,0),\,z\geq 0}\max_{t\geq0}I_{\lambda_n,\mu_n}(tz)
\vspace{0.2cm}
\\
&\leq& \dis\inf_{z \neq (0,0),\,z\geq 0}\max_{t \geq0}I_{0,0}(tz) = c_{0,0},
\end{array}
$$
where we have used, in the last equality, that the infimum $c_{0,0}$ is attained at a nonnegative solution. The above inequality  implies that
\begin{equation} \label{minimax1}
\limsup_{n\to\infty} c_{\lambda_n,\mu_n} \leq c_{0,0}.
\end{equation}

On the other hand, it follows from Theorem \ref{th_existencia} that there exists $(z_n)=(u_n,v_n) \subset X$ such that
$$
I_{\lambda_n,\mu_n}(z_n) = c_{\lambda_n,\mu_n},~~~I_{\lambda_n,\mu_n}'(z_n)=0.
$$
Since $c_{\lambda_n,\mu_n}$ is bounded, the same argument performed in (\ref{psbound}) implies that $(z_n)$ is bounded in $X$. Thus
$\int Q_{\lambda_n,\mu_n}(z_n) \leq \lambda_n \int (|u_n|^q+|v_n|^q)$, from which it follows that
\begin{equation} \label{minimax2}
\lim_{n\to\infty} \int Q_{\lambda_n,\mu_n}(z_n) = 0.
\end{equation}

Let $t_n>0$ be such that $t_n z_n \in \mathcal{N}_{0,0}$. Since $z_n \in \mathcal{N}_{\lambda_n,\mu_n}$, we have that
$$
\begin{array}{lcl}
c_{0,0}&\leq& I_{0,0}(t_nz_n)=I_{\lambda_n,\mu_n}(t_nz_n)+t_n^q\dis\int Q_{\lambda_n,\mu_n}(z_n) \vdois \\
 &\leq&  I_{\lambda_n,\mu_n}(t_nz_n)+t_n^q\dis\int Q_{\lambda_n,\mu_n}(z_n) \vdois \\
& = & c_{\lambda_n,\mu_n}+t_n^q\dis\int Q_{\lambda_n,\mu_n}(z_n).
\end{array}
$$
If $(t_n)$ is bounded, we can use the above
estimate and (\ref{minimax2}) to get
$$
c_{0,0} \leq \liminf_{n\to\infty} c_{\lambda_n,\mu_n}.
$$
This and (\ref{minimax1}) proves the lemma.

It remains to check that $(t_n)$ is bounded. A straightforward calculation shows that
\begin{equation} \label{minimax3}
t_n = \left( \frac{\|z_n\|^2}{\int H(z_n)} \right)^{1/(2^*-2)}.
\end{equation}
Since $z_n \in \mathcal{N}_{\lambda_n,\mu_n}$ we obtain
$$
\| z_n\|^2 = q \int Q_{\lambda_n,\mu_n}(z_n) + \int H(z_n) \leq o_n(1) +
S_H^{-2^*/2} \|z_n\|^{2^*}.
$$
Hence $\|z_n\|^2 \geq c_1>0$, and therefore it follows from the
 above expression that $\int H(z_n) \geq c_2>0$. This, the
boundedness of $(z_n)$ and  (\ref{minimax3}) imply that $(t_n)$ is
bounded. The lemma is proved. \qed

\section{Proof of the main theorems}

From now on we fix $r>0$  such  that the sets
$$
\Omega_r^+ := \{ x \in \mathbb{R}^N :
\mbox{dist}(x,\Omega) < r\},~~~\Omega_r^- := \{x \in
\Omega : \mbox{dist}(x,\partial\Omega) > r\}
$$
are homotopic equivalents to $\Omega$. We define the functional
 $$
 J_{\lambda,\mu} := \frac{1}{2} \| z\|^2 - \int Q_{\lambda,\mu}(z) - \frac{1}{2^*} \int H(z),~~~z \in X_{r,rad},
 $$
where $X_{r,rad} := \{ (u,v) :  u,\,v \in H_{0}^1(B_r(0)) \mbox{ and } u,\, v \mbox{ are radial functions} \}$.

We denote by $\mathcal{M}_{\lambda,\mu}$ its
associated Nehari manifold and set
$$
m_{\lambda,\mu} := \inf_{z \in \mathcal{M}_{\lambda,\mu}}
J_{\lambda,\mu}(z).
$$
According to \cite[Lemma 1]{FilSou} the infimum $S_H$ can be attained
by functions belonging to $\mathcal{D}^{1,2}_{rad}(\mathbb{R}^N) \times
\mathcal{D}^{1,2}_{rad}(\mathbb{R}^N)$. So, arguing as in the proof of Lemma \ref{lema-limite-zero} and
Theorems \ref{th_existencia} and \ref{th_existencia2} , we
obtain the following result.

\begin{lemma} \label{th_existencia-radial}
Under the hypothesis of Theorem \ref{th1}, the infimum $m_{\lambda,\mu}$ is attained by a positive radial
function $z_{\lambda,\mu} \in X_{r,rad}$ whenever $2<q<2^*$ and
$\lambda,\,\mu>0$, or $q=2$ and $\lambda,\,\mu \in
(0,\frac{\lambda_1}{2})$. Moreover
$$
m_{\lambda,\mu} < \frac{1}{N} S_H^{N/2}~~\mbox{and}~~\lim_{\lambda,\,\mu \to 0^+} m_{\lambda,\mu} = \frac{1}{N} S_H^{N/2}.
$$
The same result hold if we assume the hypothesis of Theorem \ref{th2} and $\lambda,\,\mu>0$. Moreover, both results hold if $N=3$ and $4<q<6$.

\end{lemma}

 We  introduce the barycenter map $\beta_{\lambda,\mu} : \mathcal{N}_{\lambda,\mu} \to \mathbb{R}^N$ as follows
 $$
 \beta_{\lambda,\mu}(z) := \frac{1}{S_H^{N/2}} \int H(z)x \,\textrm{d}x.
 $$
This maps has the following property.

\begin{lemma} \label{lema-baricentro}
There exists $\lambda^*>0$ such that $\beta_{\lambda,\mu}(z) \in \Omega_{r/2}^+$ whenever $z \in \mathcal{N}_{\lambda,\mu}$, $\lambda,\,\mu \in (0,\lambda^*)$ and $I_{\lambda,\mu}(z) \leq m_{\lambda,\mu}$.
\end{lemma}

\proof Suppose, by contradiction, that there exist $(\lambda_n),\,(\mu_n) \subset \mathbb{R}^+$ and $(w_n) \subset \mathcal{N}_{\lambda_n,\mu_n}$ such that $\lambda_n,\,\mu_n \to 0^+$ as $n\to\infty$, $I_{\lambda_n,\mu_n}(w_n) \leq m_{\lambda_n,\mu_n}$ but $\beta_{\lambda_m,\mu_n}(w_n) \not\in \Omega_{r/2}^+$.

Standard calculations show that $(w_n)$ is bounded in $X$. Moreover
$$
0 = I_{\lambda_n,\mu_n}'(w_n)w_n = \|w_n\|^2 - q \int Q_{\lambda_n,\mu_n}(w_n) - \int H(w_n).
$$
As in the proof of Lemma \ref{lema-limite-zero} we have that $\int
Q_{\lambda_n,\mu_n}(w_n) \to 0$ and therefore $\lim_{n\to\infty} \|w_n\|^2 =
\lim_{n\to\infty} \int H(w_n) = b \geq 0$. Notice that
$$
c_{\lambda_n,\mu_n} \leq I_{\lambda_n,\mu_n}(w_n) = \frac{1}{2}
\|w_n\|^2 - \int Q_{\lambda_n,\mu_n}(w_n) - \frac{1}{2^*} \int H(w_n) \leq
m_{\lambda_n,\mu_n}.
$$
Recalling that  $c_{\lambda_n,\mu_n}$ and $m_{\lambda_n,\mu_n}$
both converge to $\frac{1}{N} S_H^{N/2}$, we can use the above expression
and $\int Q_{\lambda_n,\mu_n}(w_n) \to 0$ again to conclude that $b=S_H^{N/2}$, that
is,
\begin{equation} \label{baricentro1}
\lim_{n\to\infty} \|w_n\|^2 = S_H^{N/2} = \lim_{n\to\infty} \int H(w_n).
\end{equation}

Let $t_n := (\int H(w_n))^{-1/2^*}>0$ and notice that $z_n:=t_nw_n$ satisfies the hypotheses of Proposition \ref{lemaSH}. Thus, for some sequences $(r_n) \subset
(0,\infty)$ and $(y_n) \subset \mathbb{R}^N$ satisfying $r_n \to
0$,  $y_n \to \overline{y} \in \overline{\Omega}$ we have that
$z_n^{y_n,r_n} \to z$ in $\mathcal{D}^{1,2}(\mathbb{R}^N) \times
\mathcal{D}^{1,2}(\mathbb{R}^N)$.

The definition of $z_n$, (\ref{baricentro1}), the strong
convergence of $(z_{n}^{y_n,r_n})$ and the Lebesgue's Theorem
provide
$$
\begin{array}{lcl}
\beta_{\lambda_n,\mu_n}(w_n) & = & \dis\frac{t_n^{-2^*}}{S_H^{N/2}} \int H(z_n)x \,\textrm{d}x = (1+o_n(1))\dis\int H(z_n)x \,\textrm{d}x \vdois \\
& = & (1+o_n(1)) \dis\int H(z_n^{y_n,r_n})(r_nx+y_n)\,\textrm{d}x \vdois\\
& = & (1+o_n(1)) \left( \dis\int H(z)\overline{y} \,\textrm{d}x + o_n(1) \right).
\end{array}
$$
Since $\overline{y} \in \overline{\Omega}$ and $\int H(z)=1$, the above expression implies that
$$
\lim_{n\to\infty}
\mbox{dist}(\beta_{\lambda_n,\mu_n}(w_n),\overline{\Omega}) = 0,
$$
which contradicts $\beta_{\lambda_n,\mu_n}(w_n) \not\in
\Omega_{r/2}^+$. The lemma is proved.\qed

According to Lemma \ref{th_existencia-radial}, for each $\lambda,\, \mu>0$ small the infimun $m_{\lambda,\mu}$ is attained by a nonnegative radial function $z_{\lambda,\mu}$. We consider
$$
I_{\lambda,\mu}^{m_{\lambda,\mu}} := \{ z \in X :
I_{\lambda,\mu}(z) \leq m_{\lambda,\mu} \}
$$
and define the function $\gamma_{\lambda,\mu} : \Omega_r^- \to I_{\lambda,\mu}^{m_{\lambda,\mu}}$ by setting, for each $y \in \Omega_r^-$,
$$
\gamma_{\lambda,\mu}(y)(x) := \left\{
\begin{array}{ll}
z_{\lambda,\mu}(x-y) & \mbox{ if } x \in B_r(y), \vdois \\
0 & \mbox{ otherwise.}
\end{array}
\right.
$$
A change of variables and straightforward calculations show that the map  $\gamma_{\lambda,\mu}$ is well defined. Since $z_{\lambda,\mu}$ is radial, we have that $\int_{B_r(0)} H(z_{\lambda,\mu})x\,\textrm{d}x=0$. Hence, for each $y \in \Omega_r^-$, we obtain
$$
\beta_{\lambda,\mu}(\gamma_{\lambda,\mu}(y)) = \alpha(\lambda,\mu)y,
$$
where
$$
\alpha(\lambda,\mu) := \frac{1}{S_H^{N/2}} \int H(z_{\lambda,\mu}).
$$

If we define $F_{\lambda,\mu}: [0,1] \times
(\mathcal{N}_{\lambda,\mu} \cap I_{\lambda,\mu}^{m_{\lambda,\mu}})
\to \mathbb{R}^N$ by
$$
F_{\lambda,\mu}(t,z) := \left( t +
\frac{1-t}{\alpha(\lambda,\mu)}\right)\beta_{\lambda,\mu}(z),
$$
we have the following.

\begin{lemma} There exists $\lambda^{**}>0$ such that,
$$
F_{\lambda,\mu}\left([0,1]\times\left(\mathcal{N}_{\lambda,\mu} \cap I^{m_{\lambda,\mu}}_{\lambda,\mu}\right)\right)\subset\Omega^+_r,
$$
whenever $\lambda,\mu\in(0,\lambda^{**})$.
\label{lemaH}
\end{lemma}

\proof Arguing by contradiction, we suppose that there exist
sequences $(\lambda_n),\,(\mu_n) \subset \mathbb{R}^+$ and
$(t_n,z_n) \in [0,1] \times (\mathcal{N}_{\lambda_n,\mu_n} \cap
I_{\lambda_n,\mu_n}^{m_{\lambda_n,\mu_n}})$  such that
$\lambda_n,\,\mu_n \to 0^+$, as $n\to\infty$, and
$F_{\lambda_n,\mu_n}(t_n,z_n) \not \in \Omega_r^+$. Up to a
subsequence $t_n \to t_0 \in [0,1]$. Moreover, the compactness of
$\overline{\Omega}$ and Lemma \ref{lema-baricentro} imply that, up
to a subsequence, $\beta_{\lambda_n,\mu_n}(z_n) \to y \in
\overline{\Omega_{r/2}^+} \subset \Omega_r^+$. We claim that
$\alpha(\lambda_n,\mu_n) \to 1$. If this is true,  we can use
the definition of $F$ to conclude that
$F_{\lambda_n,\mu_n}(t_n,z_n) \to y \in \Omega_r^+$, which does
not make sense.

It remains to check the above claim. It follows from Lemma
\ref{th_existencia-radial} that
$$
m_{\lambda_n,\mu_n} =  \frac{1}{2} \| z_{\lambda_n,\mu_n}\|^2 -
\int_{B_r(0)} Q_{\lambda_n,\mu_n}(z_{\lambda_n,\mu_n})- \frac{1}{2^*} \int_{B_r(0)}
H(z_{\lambda_n,\mu_n}) < \frac{1}{N} S_H^{N/2}.
$$
As before $\int_{B_r(0)} Q_{\lambda_n,\mu_n}(z_{\lambda_n,\mu_n}) \to 0$. This,
$J_{\lambda_n,\mu_n}'(z_{\lambda_n,\mu_n})=0$, the above
expression and  the same arguments used in the proof of Lemma \ref{th_existencia-radial}  imply that
$$
\lim_{n\to \infty} \int H(z_{\lambda_n,\mu_n}) = S_H^{N/2}.
$$
The  equality above and the definition of $\alpha(\lambda,\mu)$ imply that $\alpha(\lambda_n,\mu_n) \to 1$. The lemma is proved. \qed

\begin{corollary} \label{cor}
Let $\Lambda := \min\{\lambda^*,\lambda^{**}\} >0$, with $\lambda^*$ and $\lambda^{**}$ given by Lemmas \ref{lema-baricentro} and \ref{lemaH}, respectively. If $Q$ is such that $\lambda,\mu \in (0,\Lambda)$ then
$$
\emph{cat}_{I_{\lambda,\mu}^{m_{\lambda,\mu}}}(I_{\lambda,\mu}^{m_{\lambda,\mu}}) \geq \emph{cat}_{\Omega}(\Omega).
$$
\end{corollary}

\proof It suffices to use Lemmas \ref{lema-baricentro} and
\ref{lemaH} and argue as in \cite[Lemma 4.3]{AlvDin}.  We omit the
details. \qed

We are now ready to prove our main results.

\vspace{0.2cm} \noindent \textit{Proof of Theorems \ref{th1} and \ref{th2}.} Let
$\Lambda >0$ be given by Corollary \ref{cor} and suppose that $Q$
is such that $\lambda,\,\mu \in (0,\Lambda)$. Using Lemma
\ref{lemaPS}  and arguing as in \cite[Lemma 4.2]{AlvDin} we can
prove that the functional $I_{\lambda,\mu}$ restricted to
$\mathcal{N}_{\lambda,\mu}$ satisfies the $(PS)_c$ condition for
all $c < \frac{1}{N}S_{H}^{N/2}$. Since $m_{\lambda,\mu}<\frac{1}{N}S_H^{N/2}$,
standard Ljusternik-Schnirelmann theory provides
cat$_{I_{\lambda,\mu}^{m_{\lambda,\mu}}}(I_{\lambda,\mu}^{m_{\lambda,\mu}})
$ critical points of the constrained functional. If $z
\in\mathcal{N}_{\lambda,\mu}$ is one of these critical points, the
same argument of \cite[Lemma 4.1]{AlvDin} shows that $z$ is also a
critical point of the unconstrained functional, and therefore a
nontrivial solution of $(P)$. As before, the obtained solutions
are nonnegative in $\Omega$.  The results follow from Corollary
\ref{cor}. \qed

\section{Some further remarks} \label{secao_final}

We start this last section presenting some functions which satisfy our hypotheses. We have the
following example from \cite{FilSou}. Let $2 \leq q < 2^*$ and
$$
P_q(s,t) := a_1 s^{q} + a_2 t^{q} + \sum_{i=1}^k b_i s^{\alpha_i}t^{\beta_i},~~~s,\,t \geq 0,
$$
where $\alpha_i, \beta_i > 1$, $\alpha_i+\beta_i=q$ and $a_1,\,a_2,\,b_i \in \mathbb{R}$. The following functions and its possible
combinations, with appropriated choices of the coefficients $a_1,\,a_2,\,b_i$,
satisfy our hypotheses on $Q$
$$
Q(s,t) = P_q(s,t),~~Q(s,t) =
\sqrt[r]{P_{rq}(s,t)}~~\mbox{and}~~Q(s,t) =
\frac{P_{r+l}(s,t)}{P_{l}(s,t)},
$$
with $l>0$. Hence, we see that our subcritical term is more general than those of \cite{Han,Jap1,Jap2}.

The form of $H$ is more restricted due to $(H_4)$. This technical condition has already appeared in \cite{FilSou} and it is important to guarantee that the constant $S_H$ defined in (\ref{def_sh}) does not depend on $\Omega$.
As quoted in \cite{FilSou}, the concavity condition $(H_4)$ is satisfied if $H \in C^2(\mathbb{R}^2_+,\mathbb{R})$ is such that $H_{st}(s,t) \geq 0$ for each $(s,t) \in \mathbb{R}^2_+$.

Although we have more restrictions on the shape of $H$, it can have the
polynomial form
$$
H(s,t)=P_{2^*}(s,t).
$$
Thus, differently from \cite{Han,Jap1,Jap2}, we can deal here with functions $H$ which possesses coupled and no coupled terms. For example, the function
$$
H(s,t)=a_1s^{2^*} + a_2t^{2^*} + a_3s^{\alpha}t^{\beta},
$$
with $a_i \in \mathbb{R}$, $\alpha,\beta > 1$, $\alpha + \beta = 2^*$ satisfies the hypotheses $(H_0)-(H_4)$ for appropriated choices of the coefficients $a_i$. We also mention that the positivity condition in $(H_2)$ can holds even if some of the coefficients $a_i$ are negative. As a simple example, suppose that $H$ is as above with $a_1,\,a_2 \geq 0$ and $a_3<0$. Since  $s^{\alpha}v^{\beta} \leq  s^{2^*}+t^{2^*}$, the condition $(H_2)$ holds for  $a_3 > \max \{ -a_1,-a_2\}$.

Another interesting remark is that we can obtain  versions of our theorems by interchanging conditions like $(Q_1)$ and $(\widehat{Q_1})$ for both the functions $Q$ and $H$. More specifically, let us consider the following assumption
\begin{itemize}
\item[($\widehat{H_1})$] $H_s(0,1)>0$ and $H_t(1,0)>0$.
\end{itemize}
A simple inspection of our proofs shows that Theorem \ref{th1} is valid if we suppose $(\widehat{H_1})$ and $(Q_1)$. The same is true for Theorem \ref{th2}. This last theorem is also true if we suppose $(\widehat{H_1})$ and $(\widehat{Q_1})$. The difference among these various settings relies in the form of the possible coupled terms.

A simple inspection of our proofs show that, instead of just one subcritical term, we can consider in $(P)$ a subcritical nonlinear term of the form
$$
\widetilde{Q}(s,t) = \sum_{i=1}^k Q_i(s,t),
$$
with each function $Q_i$ being $q_i$-homogeneous, $2 \leq q_i < 2^*$, and satisfying the same kind of hypotheses of $Q$. In this case, for each $i=1,\ldots,k$, we define the numbers $\mu_i,\,\lambda_i$ as in (\ref{defmu})-(\ref{deflambda}), and the results hold if $\max_{i=1,\ldots,k}\{\mu_i,\lambda_i\}$ is small enough.

With some additional conditions we can assure that the solutions obtained in this paper are positive. Indeed, if we suppose that
\begin{itemize}
\item[$(Q_2)$] $Q_s(s,t) \geq0,\,Q_t(s,t) \geq 0$ for each $(s,t) \in \mathbb{R}^2_+$,
\end{itemize}
we can apply the Maximum Principle in each equation of $(P)$. Thus, if $(u,v)$ is a nonnegative solution, then $u \equiv 0$ or $u>0$ in $\Omega$, the same holding for $v$. We need only to discard solutions of the type $(u,0)$ or $(0,v)$. This can be done if we guarantee some kind of strongly coupling for the system. In what follows, we present some situations where this can be done.

If we are under the conditions of Theorem \ref{th1} we assume a stronger form of $(Q_1)$ and $(H_1)$, namely that $\nabla Q(1,0)=\nabla Q(0,1) = \nabla H(1,0) = \nabla H(0,1)=0$. In this way, if $(u,0)$ is a solution then
$$
0=I_{\lambda,\mu}'(u,0)(u,0) = - \| u\|^2 - \int \left( Q_u(u,0)u + \frac{1}{2^*} H_u(u,0)u \right) = - \|u\|^2.
$$
and therefore $u \equiv 0$. Analogously, if $(0,v)$ is a solution then $v \equiv 0$. In the setting of Theorem \ref{th2} and considering the solution $(u,0)$ we obtain, from the second equation, that
$$
0 = Q_v(u,0) + H_v(u,0) = u^{q-1} Q_v(1,0).
$$
Since from $(\widehat{Q_1})$ we have that $Q_v(1,0)>0$, it follows that $u \equiv 0$. The argument for $(0,v)$ is analogous.

\end{document}